# Asymptotics of the solution of the turbulent diffusion equation taking into account the polydispersity of the impurity and wind pickup from the underlying surface[*]

**Nesterov A.V.** [1[0000-0002-4702-4777]]

[1] PLEKHANOV Russian University of Economics, Stremyanny lane 36, Moscow, 117997, Russia
andrenesterov@yandex.ru
[2]

**Abstract.** The asymptotics of a singularly perturbed problem is constructed. describing the transport of a polydisperse impurity in the atmosphere, taking into account the processes of precipitation and wind pick-up, as well as the processes of coagulation - dissociation. The mathematical model of this process represents a differential-operator equation of turbulent diffusion with a non-standard boundary condition containing two components - atmospheric and soil. The asymptotics of the solution is constructed by the method of boundary functions. Problems that do not contain small parameters are obtained for the main terms of the asymptotic equation.

## 1 Introduction

We develop a model of impurity transport in the atmosphere, taking into account the following factors.

1. The impurity is carried by the wind flow in the horizontal direction and diffuses in the horizontal direction.

2. The impurity settles and diffuses in the vertical direction.

3. The impurity settles on the underlying surface and is picked up from the underlying surface.

4. The impurity is polydisperse, and both coagulation and dissociation of the impurity in the air occur.

---

[*] This research was performed in the framework of the state task in the field of scientific activity of the Ministry of Science and Higher Education of the Russian Federation, project "Models, methods, and algorithms of artificial intelligence in the problems of economics for the analysis and style transfer of multidimensional datasets, time series forecasting, and recommendation systems design", grant no. FSSW-2023-0004. design", grant no. FSSW-2023-0004.

The figure shows an illustration of the processes taking place.

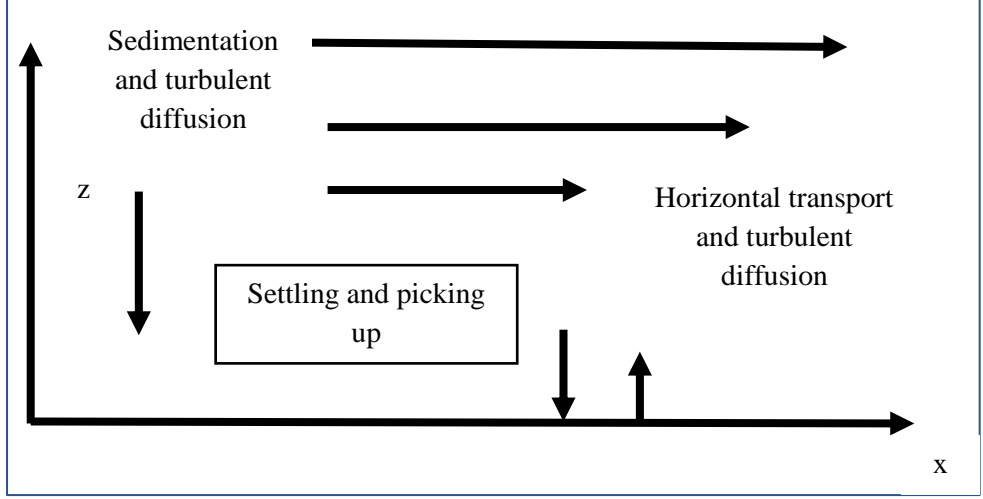

As is known, impurity transport in the atmospheric boundary layer is described quite well by the turbulent diffusion equation, taking into account the dependence of the transport rates and diffusion coefficients on the height above the underlying surface [1]. Directing the x-coordinate along the wind direction, and the z - coordinate vertically upwards, taking into account the change in the transport coefficients with altitude, neglecting changes along the horizontal and considering the transport process stationary, we write this equation in the form:

$$u_t + V(z)\frac{\partial u}{\partial x} - w\frac{\partial u}{\partial z} = K_x(z)\frac{\partial^2 u}{\partial x^2} + \frac{\partial}{\partial z}\left(K_z(z)\frac{\partial u}{\partial z}\right),$$

$$|x| < \infty,\ z > 0, t > 0,$$

(1)

where u(x,z,t) is the impurity concentration in the atmosphere, V(z) is the wind speed, w is the deposition speed, $K_x(z)$ and $K_z(z)$ are the coefficients of turbulent diffusion in the horizontal and vertical directions, respectively.

Equation (1) must be supplemented with the initial and boundary conditions. The initial condition is set in the usual way, and the boundary condition in the absence of subsidence and lifting from the underlying surface is usually set in the form of the Monin condition [2]

$$wu + K_z(0)\frac{\partial u}{\partial z}\bigg|_{z=0} = \alpha u\big|_{z=0},$$

(2)

where the coefficient α is the so-called dry deposition rate. The Monin boundary condition has been studied in many papers, in particular, in [3].

Condition (2) does not describe the processes of impurity lifting from the surface, so in the case of a "dusty surface", it should be replaced with a more complex one that takes into account the impurity lying on the surface. Such conditions were proposed, for example, in [4], [5].

Following [4], [5], we describe the process of deposition-wind pickup of systemsй of linear equations-of boundary conditions

$$v_t = \left( wu + K_z(0)\frac{\partial u}{\partial z} \right)\bigg|_{z=0},$$

$$v_t = \alpha u|_{z=0} - \beta v,$$

<div align="right">(3)</div>

where v(x,t)is the impurity concentration on the soil surface, and α and β are the so-called коэффициенты dry deposition and pick-up coefficients, respectively. The physical meaning of these coefficients is described in detail, for example ,in [5].

The first relation in (3) gives the mass balance in the soil - atmosphere system ( the rate of change in the impurity concentration on the soil surface is equalto the flow from the atmosphere). The second relation in (3) expressesь the rate of change in the concentration on the soil surface as разностьthe difference between two fluxes - precipitation of an impurity from the atmosphere $\alpha u|_{z=0}$ and discharge of an impurity into the atmosphere βv ( fluxes are considered linearly dependent on the concentrationof x). For β=0, the system of conditions passes into the well-known Monin condition (2). Previously, these conditions for β≠0 were used, for example, in работах [6], [7].

In [6], [7], the dispersed composition of the impurity was not taken into account. Taking into account the impurity dispersion leads to the dependence of the parameters of the transfer process on the size of the particles containing the impurity.

We assume that the impurity in the atmosphere is polydisperse, i.e., it consists of particles of various sizes. In this case, it is advisable to add an additional variable to the variables that depend on the impurity concentration, and consider u=u(x,z,t,p), where p is the size (mass) particles, p⊂P, where P is the set of acceptable sizes (masses), for example,P=[$p_{p1}$,$p_{p2}$].

Let coagulation-dissociation processes occur in the atmosphere, which are described by the equation

$$\frac{\partial u(t,p)}{\partial t} = L_p u$$

<div align="right">(4)</div>

where $L_p$ is a linear operatoracting on the variable p. If there is an equilibrium size distribution of the impurity, $h_0$h0(t,p), then there is a distribution for which

$$\frac{\partial h_0(t,p)}{\partial t} \equiv 0,$$

<div align="right">(5)</div>

this means that the operator $L_p$ has a zero eigenvalue, which corresponds to the eigenfunction $h_0$(p).

Remark. The $L_p$ operator can be, for example, integral

$$L_p u = \int_P K(p,q)u(q)dq$$

In this paper, it is shown that coagulation-dissociation processes occur only in the atmosphere.

We assume that the horizontal transport and diffusion coefficients, as well as the deposition-uplift parameters, depend on the particle size p. In this case, the variables u,v will depend on one more variable p, u=u(x,z,t,p), v=v(x,t,p), part of the transfer coefficients will also depend on the parameter $p$

$$u_t + V(z)\frac{\partial u}{\partial x} - w(p)\frac{\partial u}{\partial z} = K_x(z)\frac{\partial^2 u}{\partial x^2} + \frac{\partial}{\partial z}\left(K_z(z,p)\frac{\partial u}{\partial z}\right) + L_p u,$$

$$|x| < \infty, z > 0, t > 0.$$

(6)

The coefficients in the boundary conditions will also be considered to depend on the parameter p

$$v_t = \left(w(p)u + K_z(0,p)\frac{\partial u}{\partial z}\right)\bigg|_{z=0},$$

$$v_t = \alpha(p)u\big|_{z=0} - \beta(p)v.$$

(7)

The system (6)-(7) should be supplemented with the initial conditions

$$u(x,z,0,p) = u^0(x,z,p), v(x,0,p) = v^0(x,p).$$

(8)

Solutions of the system (6)-(7) satisfy the conservation law.

$$\left(Mu(p) + Mv(p), h*(p)\right) = const =$$

$$= \left(Mu^0(x,z,p) + Mv^0(x,p), h*(p)\right), t > 0.$$

(9)

Relations (9) can be called the law of conservation of "generalized mass".

In this paper, we investigate the effect of impurity dispersion on the transfer process under the assumptions formulated below, which are quite physically justified.

The main processes described by equation (6) and conditions (7), transport and diffusion in the wind direction (along the x-axis), precipitation and diffusion along the vertical (along the z-axis), precipitation to the surface-secondary ascent, as well as coagulation-dissociation-can occur with different characteristic times, which are indicated by $T_{gor}$, $T_{vert}$, $T_{resusp}$, $T_{coagul}$.

## 2    Problem statement

The following case is considered

$T_{gor} \sim T_{vert} >> T_{resusp} \sim T_{coagul.}$

When the equation is dimensioned, it becomes singularly perturbed [8]:

$$\varepsilon^2\left(u_t + V(z)\frac{\partial u}{\partial x} - w(p)\frac{\partial u}{\partial z} - K_x(z)\frac{\partial^2 u}{\partial x^2} - \frac{\partial}{\partial z}\left(K_z(z,p)\frac{\partial u}{\partial z}\right)\right) = L_p u,$$

$$|x| < \infty, z > 0, t > 0,$$

$$\left(w(p)u + K_z(0,p)\frac{\partial u}{\partial z}\right)\bigg|_{z=0} = v_t,$$

$$\varepsilon^2 v_t = \alpha(p)u\big|_{z=0} - \beta(p)v,$$

$$u(x,z,0,p) = u^0(x,z,p), v(x,0,p) = v^0(x,p).$$

(10a,b,c,d,e)

Here $\varepsilon = T_{resusp}/T_{vert} << 1$ is a small positive parameter.

Equation (10a) refers to singularly perturbed equations in the so-called critical case [8], i.e., when the degenerate equation (for $\varepsilon = 0$)

$$L_p u = 0$$

it has a family of solutions

$$u = C(x,z,t)h_0(p).$$

Note that a similar problem for a singularly perturbed differential operator equation was investigated, for example ,in [10], but with a different arrangement of small parameters.

The main goal of this paper is to construct the principal terms of the formal asymptotic expansion (PAR) of the solution of a singularly perturbed problem and (10a,b,c,d,e) with respect to a small parameter $\varepsilon$

$$\varepsilon^2\left(u_t + V(z)\frac{\partial u}{\partial x} - w(p)\frac{\partial u}{\partial z} - K_x(z)\frac{\partial^2 u}{\partial x^2} - \frac{\partial}{\partial z}\left(K_z(z,p)\frac{\partial u}{\partial z}\right)\right) = L_p u,$$

$$|x| < \infty, z > 0, t > 0,$$

$$\left(w(p)u + K_z(0,p)\frac{\partial u}{\partial z}\right)\bigg|_{z=0} = v_t,$$

$$\varepsilon^2 v_t = \alpha(p)u|_{z=0} - \beta(p)v,$$

$$u(x,z,0,p) = u^0(x,z,p), v(x,0,p) = v^0(x,p).$$

## 3  Construction of an AE

We restrict ourselves to constructing the highest terms of the asymptotic function with respect to a small parameter (the first non-zero terms).

Below, we use an algorithm for constructing PAR solutions to singularly perturbed problems in the critical case [8].

### 3.1. The regular part of the AR solution $\overline{u}(x,z,t,p), \overline{v}(x,t,p)$.

We will search for the regular part of the PAR of the solution of problem (10) in the standard form [8]

$$\overline{u}(x,z,t,p) = \overline{u}_0(x,z,t,p) + \varepsilon^2\overline{u}_2(x,z,t,p) + ...,$$

$$\overline{v}(x,t,p) = \overline{v}_0(x,t,p) + \varepsilon^2\overline{v}_2(x,t,p) + ....$$

Substituting the expansion for the function $\overline{u}(x,z,t,p)$ in equation (10a) and separating the powers of the parameter $\varepsilon$, we obtain:

By

$$\varepsilon^0 : 0 = L_p\overline{u}_0,$$

$$\left(w(p)\overline{u}_0 + K_z(0,p)\frac{\partial\overline{u}_0}{\partial z}\right)\bigg|_{z=0} = \overline{v}_{0,t},$$

$$0 = \alpha(p)\overline{u}_0|_{z=0} - \beta(p)\overline{v}_0.$$

<div align="right">(11)</div>

by

$$\varepsilon^2 : L_p \overline{u}_2 = l_2 =$$

$$= \left( \overline{u}_{0,t} + V(z)\frac{\partial \overline{u}_0}{\partial x} - w(p)\frac{\partial \overline{u}_0}{\partial z} - K_x(z)\frac{\partial^2 \overline{u}_0}{\partial x^2} - \frac{\partial}{\partial z}\left( K_z(z,p)\frac{\partial \overline{u}_0}{\partial z} \right) \right),$$

$$\left( w(p)\overline{u}_2 + K_z(0,p)\frac{\partial \overline{u}_2}{\partial z} \right)\Bigg|_{z=0} = \overline{v}_{2,t},$$

$$\overline{v}_{0,t} = \alpha(p)\overline{u}_2\big|_{z=0} - \beta(p)\overline{v}_2.$$

$$(12)$$

By

$$\varepsilon^k : L_p \overline{u}_k = l_2 =$$

$$= \left( \overline{u}_{k-2,t} + V(z)\frac{\partial \overline{u}_{k-2}}{\partial x} - w(p)\frac{\partial \overline{u}_{k-2}}{\partial z} - K_x(z)\frac{\partial^2 \overline{u}_{k-2}}{\partial x^2} - \frac{\partial}{\partial z}\left( K_z(z,p)\frac{\partial \overline{u}_{k-2}}{\partial z} \right) \right),$$

$$\left( w(p)\overline{u}_k + K_z(0,p)\frac{\partial \overline{u}_k}{\partial z} \right)\Bigg|_{z=0} = \overline{v}_{k,t},$$

$$\overline{v}_{k-2,t} = \alpha(p)\overline{u}_k\big|_{z=0} - \beta(p)\overline{v}_k.$$

From the first equation (11) we obtain

$$\overline{u}_0(x,z,t,p) = h_0(p)\varphi_0(x,z,t),$$

$$(13)$$

where $\varphi_0(x,z,t)$ is an unknown function yet.

Let us write down the solvability conditions of the first equation in (12)

$$(l_1, h^*_0) = 0 =$$

$$= \left( h_0\varphi_{0,t} + V(z)\frac{\partial h_0\varphi_0}{\partial x} - w(p)\frac{\partial h_0\varphi_0}{\partial z} - K_x(z)\frac{\partial^2 h_0\varphi_0}{\partial x^2} - \frac{\partial}{\partial z}\left( K_z(z,p)\frac{\partial h_0\varphi_0}{\partial z} \right), h^*_0 \right) =$$

$$= (h_0, h^*_0)\varphi_{0,t} + V(z)(h_0, h^*_0)\frac{\partial \varphi_0}{\partial x} - (h_0 w, h^*_0)\frac{\partial \varphi_0}{\partial z} - K_x(z)(h_0, h^*_0)\frac{\partial^2 \varphi_0}{\partial x^2} -$$

$$- \frac{\partial}{\partial z}\left( (h_0 K_z(z,p), h^*_0)\frac{\partial \varphi_0}{\partial z} \right),$$

From here

$$\varphi_{0,t} + V(z)\frac{\partial \varphi_0}{\partial x} - w_{ef}\frac{\partial \varphi_0}{\partial z} - K_x(z)\frac{\partial^2 \varphi_0}{\partial x^2} - \frac{\partial}{\partial z}\left( K_{z,ef}(z)\frac{\partial \varphi_0}{\partial z} \right) = 0,$$

$$w_{ef} = (h_0 w, h^*_0), K_{z,ef}(z) = (h_0 K_z(z,p), h^*_0)$$

From the third equation (11) we obtain

$$\overline{v}_0(x,t,p) = \frac{a(p)}{b(p)}\overline{u}_0(x,t,0,p) = c(p)h_0(p)\varphi_0(x,0,t),$$

$$c(p) = \frac{a(p)}{b(p)}.$$

(13)

Substituting this expression into the second equation (11), we obtain

$$\left. \left( w(p)\overline{u}_0 + K_z(0,p)\frac{\partial \overline{u}_0}{\partial z} \right) \right|_{z=0} = c(p)\,\overline{u}_{0,t}\big|_{z=0}.$$

$$\left. \left( w(p)h_0(p)\varphi_0 + K_z(0,p)\frac{\partial h_0(p)\varphi_0}{\partial z} \right) \right|_{z=0} = c(p)h_0(p)\varphi_{0,t}\big|_{z=0}.$$

Multiplying by h*, we obtain the boundary condition for $\varphi_0$

$$\left. \left( (w(p)h_0, h_0^*)\varphi_0 + (K_z(0,p)h_0, h_0^*)\frac{\partial \varphi_0}{\partial z} \right) \right|_{z=0} = (c(p)h_0, h_0^*)\,\varphi_{0,t}\big|_{z=0}.$$

or

$$\left. \left( w_{ef}\varphi_0 + K_{z,ef}(z)\frac{\partial \varphi_0}{\partial z} \right) \right|_{z=0} = c_{ef}\,\varphi_{0,t}\big|_{z=0},$$

$$c_{ef} = (c(p)h_0, h_0^*).$$

$\varphi_0$ defined as the solution of the equation

$$\varphi_{0,t} + V(z)\frac{\partial \varphi_0}{\partial x} - w_{ef}\frac{\partial \varphi_0}{\partial z} - K_x(z)\frac{\partial^2 \varphi_0}{\partial x^2} - \frac{\partial}{\partial z}\left( K_{z,ef}(z)\frac{\partial \varphi_0}{\partial z} \right) = 0,$$

$$w_{ef} = (h_0 w, h_0^*), K_{z,ef}(z) = (h_0 K_z(z,p), h_0^*)$$

with the boundary condition

$$\left. \left( w_{ef}\varphi_0 + K_{z,ef}(z)\frac{\partial \varphi_0}{\partial z} \right) \right|_{z=0} = c_{ef}\,\varphi_{0,t}\big|_{z=0},$$

$$c_{ef} = (c(p)h_0, h_0^*).$$

The second component is determined from the relations

$$\overline{v}_0(x,t,p) = \frac{a(p)}{b(p)}\overline{u}_0(x,t,0,p) = c(p)h_0(p)\varphi_0(x,0,t),$$

$$c(p) = \frac{a(p)}{b(p)}.$$

The remaining expansion terms are constructed in the standard way and are not given here. $\overline{u}_0(x,t,0,p)$ and $\overline{v}_0(x,t,p) = c(p)\overline{u}_0(x,t,0,p)$ they can't meet the initial conditions at the same time.

### 3.2. _Borderline functions_ $\Pi u(x, z, \tau, p), \Pi v(x, \tau, p)$

To satisfy the initial conditions, we construct boundary functions

$$\Pi u(x, z, \tau, p) = \Pi_0 u(x, z, \tau, p) + \varepsilon^2 \Pi_2 u(x, z, \tau, p) + \ldots,$$

$$\Pi v(x, \tau, p) = \Pi_0 v(x, \tau, p) + \varepsilon^2 \Pi_2 v(x, \tau, p) + \ldots,$$

$$\tau = \frac{t}{\varepsilon^2}. \tag{14}$$

The boundary functions $\Pi u$, $\Pi v$ must, together with the regular part of the solution, satisfy the initial conditions

$$\bar{u}(x, z, t, p) + \Pi u(x, z, \tau, p)\big|_{t=0} = u^0(x, z, p),$$

$$\bar{v}(x, t, p) + \Pi v(x, \tau, p)\big|_{t=0} = v^0(x, p), \tag{15}$$

be functions of the boundary layer, i.e. tend to zero as $\tau$ tends to plus infinity

$$\Pi u(x, z, \tau, p) \underset{\tau \to +\infty}{\longrightarrow} 0,$$

$$\Pi v(x, \tau, p) \underset{\tau \to +\infty}{\longrightarrow} 0, \tag{16}$$

satisfy the equations or parts of

$$\Pi u_\tau + \varepsilon^2 \left( V(z) \frac{\partial \Pi u}{\partial x} - w(p) \frac{\partial \Pi u}{\partial z} - K_x(z) \frac{\partial^2 \Pi u}{\partial x^2} - \frac{\partial}{\partial z}\left( K_z(z, p) \frac{\partial \Pi u}{\partial z} \right) \right) =$$

$$= L_p \Pi u,$$

$$|x| < \infty, z > 0, \tau > 0,$$

$$\left( w(p) \Pi u + K_z(0, p) \frac{\partial \Pi u}{\partial z} \right)\Big|_{z=0} = \varepsilon^{-2} \Pi v_\tau,$$

$$\Pi v_\tau = \alpha(p) \Pi u\big|_{z=0} - \beta(p) \Pi v,$$

The second and third equations cannot be fulfilled simultaneously, so we discard the boundary condition with respect to z.

Substituting expansions of border functions in the equations

$$\Pi u_\tau + \varepsilon^2 \left( V(z) \frac{\partial \Pi u}{\partial x} - w(p) \frac{\partial \Pi u}{\partial z} - K_x(z) \frac{\partial^2 \Pi u}{\partial x^2} - \frac{\partial}{\partial z}\left( K_z(z, p) \frac{\partial \Pi u}{\partial z} \right) \right) =$$

$$= L_p \Pi u,$$

$$|x| < \infty, z > 0, \tau > 0,$$

$$\Pi v_\tau = \alpha(p) \Pi u\big|_{z=0} - \beta(p) \Pi v,$$

we get

$$\Pi u_{0,\tau} = L_p \Pi u_0, \tau > 0,$$

$$\Pi v_{0,\tau} = \alpha(p) \Pi u_0\big|_{z=0} - \beta(p) \Pi v_0,$$

for even k's

$$\Pi u_{k,\tau} = L_p \Pi u_k + \pi u_k, \tau > 0,$$

$$\Pi v_{k,\tau} = \alpha(p)\Pi u_k\big|_{z=0} - \beta(p)\Pi v_k + \pi v_k,$$

where $\pi u_k, \pi v_k$ are expressed using $\Pi u_{j,\tau}$, $\Pi v_{j,\tau}$, $j < k$.

For odd k's $\Pi u_k = 0$.

Impose condition I.

*Condition 1.* Будем We assume that the operator $L_p$ has a simple discrete spectrum $\lambda_k$, $k=0,1,...$; $\lambda_0 = 0$, $Re$ $\lambda_k < 0$, $k=1,2,...$, $\lambda_k \leftrightarrow h_k(p)$.

Solving the equation

$$\Pi u_{0,\tau} = L_p \Pi u_0, \tau > 0,$$

we get

$$\Pi u_{0,\tau} = \sum_{i=0}^{\infty} g_i(x,z) h_i(p) e^{\lambda_i t}.$$

Based on

$$\Pi u(x,z,\tau,p) \underset{\tau \to +\infty}{\to} 0$$

we get

$$\Pi u_{0,\tau} = \sum_{i=1}^{\infty} g_i(x,z) h_i(p) e^{\lambda_i t}.$$

Decompose the initial condition into a series of eigenfunctions

$$u^0(x,z,p) = \sum_{i=0}^{\infty} u^0_i(x,z) h_i(p)$$
.

Substituting the decomposition into the initial condition and decomposing

$$\varphi_0(x,z,0)h_0(p) + \Pi u_0(x,z,0,p) = u^0(x,z,p) = \sum_{i=0}^{\infty} u^0_i(x,z)h_i(p)$$

$$c(p)\varphi_0(x,0,0)h_0(p) + \Pi v_0(x,0,p) = v^0(x,p),$$

$$\varphi_0(x,z,0)h_0(p) + \sum_{i=1}^{\infty} g_i(x,z)h_i(p) = \sum_{i=0}^{\infty} u^0_i(x,z)h_i(p),$$

we get

$$\varphi_0(x,z,0) = u^0_0(x,z), g_i(x,z) = u^0_i(x,z), i = 1,2,...$$

$$c(p)u^0_0(x,0)h_0(p) + \Pi v_0(x,0,p) = v^0(x,p), \Rightarrow$$

$$\Pi v_0(x,0,p) = v^0(x,p) - c(p)u^0_0(x,0)h_0(p).$$

The initial conditions for the rest $\Pi u_i, \Pi v_i$ are constructed in the same way.

Function $\Pi u_0$ found explicitly

$$\Pi u_0 = \sum_{i=1}^{\infty} u^0_i(x,z) h_i(p) e^{\lambda_i t}.$$

$\Pi v_0$ are determined from the equation

$$\Pi v_{0,\tau} = \alpha(p)\Pi u_0\big|_{z=0} - \beta(p)\Pi v_0,$$

with the initial conditionем

$$\Pi v_0(x,0,p) = v^0(x,p) - c(p)u^0_0(x,0)h_0(p).$$

The remaining expansion terms $\Pi u$, $\Pi v$ are constructed in the standard way.

Border $\Pi u$, $\Pi v$ functions that satisfy the equations and initial conditions, but do not satisfy the boundary condition.

We will impose an additional condition on the initial data.

Condition 2.

$$\alpha(p)u^0(x,z,p)\big|_{z=0} - \beta(p)v^0(x,p) = 0.$$

Note that if condition 1 is met, functions $\Pi u_0$, $\Pi v_0$ с точностью $O(\varepsilon^2)$ that do not accurately introduce an error in the boundary condition.

Finally, the AR for solving the problem (1)-(2) with accuracy $O(\varepsilon^2)$ is constructed in the form

$$u(x,z,t,p) = \overline{u}_0(x,z,t,p) + \Pi_0 u(x,z,\tau,p) + Ru,$$

$$v(x,t,p) = \overline{v}_0(x,t,p) + \Pi_0 v(x,\tau,p) + Rv,$$

## 4        Evaluation of the residual term

The evaluation of the residual term is based on the discrepancy.

*The theorem.* If the conditions 1,2 residual terms $Ru, Rv$ are met there is a solution to the problem

$$\varepsilon^2\left(Ru_t + V(z)\frac{\partial Ru}{\partial x} - w(p)\frac{\partial Ru}{\partial z} - K_x(z)\frac{\partial^2 Ru}{\partial x^2} - \frac{\partial}{\partial z}\left(K_z(z,p)\frac{\partial Ru}{\partial z}\right)\right) =$$

$$= L_p Ru + O(\varepsilon^2),$$

$$|x| < \infty, z > 0, t > 0,$$

$$\left(w(p)Ru + K_z(0,p)\frac{\partial Ru}{\partial z}\right)\Big|_{z=0} = Rv_t + O(\varepsilon^2),$$

$$\varepsilon^2 Rv_t = \alpha(p)Ru\big|_{z=0} - \beta(p)Rv + O(\varepsilon^2),$$

$$Ru(x,z,0,p) = 0, Rv(x,0,p) = 0.$$

## 4        Results

1. An AR solution to problem (10) is constructed.

2. The main members of the AP that are most interesting from the point of view of applications $\overline{u}_0(x,z,t,p), \overline{v}_0(x,t,p)$ are defined by the formulas

$$\overline{u}_0(x,z,t,p) = h_0(p)\varphi_0(x,z,t),$$

$$\overline{v}_0(x,t,p) = c(p)h_0(p)\varphi_0(x,0,t), c(p) = \frac{a(p)}{b(p)}.$$

Here the function $\varphi_0(x,z,t)$ is problem solving

$$\varphi_{0,t} + V(z)\frac{\partial \varphi_0}{\partial x} - w_{ef}\frac{\partial \varphi_0}{\partial z} - K_x(z)\frac{\partial^2 \varphi_0}{\partial x^2} - \frac{\partial}{\partial z}\left(K_{z,ef}(z)\frac{\partial \varphi_0}{\partial z}\right) = 0,$$

$$w_{ef} = (h_0 w, h^*_0), K_{z,ef}(z) = (h_0 K_z(z,p), h^*_0),$$

$$\varphi_0(x,z,0) = u^0_0(x,z).$$

$$\left.\left(w_{ef}\varphi_0 + K_{z,ef}(z)\frac{\partial \varphi_0}{\partial z}\right)\right|_{z=0} = c_{ef}\left.\varphi_{0,t}\right|_{z=0},$$

$$c_{ef} = (c(p)h_0, h^*_0).$$

Border $\Pi_0 u(x,z,\tau,p), \Pi_0 v(x,\tau,p)$ functions are defined by

$$\Pi u_{0,\tau} = \sum_{i=1}^{\infty} g_i(x,z)h_i(p)e^{\lambda_i t}.$$

$$\Pi v_{0,\tau} = \alpha(p)\Pi u_0|_{z=0} - \beta(p)\Pi v_0,$$

$$\Pi v_0(x,0,p) = v^0(x,p) - c(p)\varphi_0(x,0,0)h_0(p),$$

where

$$u^0(x,z,p) = \sum_{i=0}^{\infty} u^0_i(x,z)h_i(p)$$

3. To use formulas, it is sufficient to know onlythe operator's own $h_0(p)$ function corresponding $L_p$ to the zero eigenvalue (equilibrium distribution), and the fulfillment of conditions 1, 2.

## References.

1. A. S. Monin, A. M. Yaglom. Statistical fluid mechanics. Turbulence mechanics. Part 1. Moscow, Nauka Publ., 1965, 640 pp. The equation

2. Monin A. S. On the boundary condition on the Earth's surface for diffusing admixture. In the collection Atmospheric Diffusion and Air Pollution, Moscow, IL., 1962

3. Byzova N. L., Krotova I. A., Natanzon G. A. On the boundary condition in problems of impurity scattering in the atmosphere.- Meteorology and Hydrology, 1980, No. 2.

4. Slinn W. Formulation and solution of the diffusion-deposition-resuspension problem. - Atmos. Environ. 1976? vol. 10 No. 3-boundary condition with precipitation

5. Transuranic elements in the environment / translated from English-Moscow, Energoatomizdat, 1985.

6. O. I. Vozzhennikov, A.V. Nesterov. On the transport of impurities in the atmosphere during wind pick-up from the underlying surface. Meteorology and Hydrology, 1988, No. 11, pp. 63-70.

7. O. I. Vozzhennikov, A.V. Nesterov. On the boundary condition for the turbulent diffusion equation for a dusty underlying surface. Meteorology and Hydrology, 1991, No. 3, pp. 32-38.

8. Vasilyeva A. B., Butuzov V. F. Singularly perturbed equations in critical cases, Moscow: Mosk Publishing House. univ., 1978, 106 p.

9. Zaborsky A.V., Nesterov A.V. Asymptotic expansion of the solution of a singularly perturbed differential operator equation in the critical case. Mathematical Modeling, 2014, vol. 26, No. 4, pp. 65-79.

10. Zaborsky A.V., Nesterov A.V. Mathematical model переноса of polydisperse impurity transport in the atmosphere during wind pickup from the underlying surface. Modern methods of the theory of functions and related problems: materials of the International Conference: Voronezh Winter Mathematical School (January 28-February 2, 2021) / VGU; MGU im. M. V. Lomonosov; МИ I im. V. A. Steklov RAS. - Voronezh: VSU Publishing House, 2021– - 333 p.